\documentclass[11pt]{amsart}
\usepackage{amsmath, amssymb, amsthm, amscd, graphicx}

\theoremstyle{plain}
\newtheorem{prop}{Proposition}[section]
\newtheorem{coro}[prop]{Corollary}
\newtheorem{conj}[prop]{Conjecture}
\newtheorem{lemm}[prop]{Lemma}
\newtheorem{ques}[prop]{Question}

\theoremstyle{definition}
\newtheorem{defi}[prop]{Definition}

\newtheorem{exam}[prop]{Example}
\newtheorem{rema}[prop]{Remark}

\numberwithin{equation}{section}

\def\Reff#1; #2; #3; #4; #5; #6; #7\par{%
\bibitem{#1} #2, {\it #3}, #4 {\bf #5} (#6) #7}

\def\Ref#1; #2; #3; #4\par{%
\bibitem{#1} #2, {\it #3}, #4}

\catcode`\¨=\active\def¨{~} 

\renewcommand{\aa}[3]{a_{#1,#2,#3}}
\newcommand{\aab}[2]{\overline{a}_{#1,#2}}
\newcommand{\aas}[3]{\widehat a_{#1, #2,#3}}

\newcommand{\bb}[2]{b_{#1,#2}}
\newcommand{\bbb}[3]{b_{#1,#2}(#3)}
\newcommand{\bbbb}[3]{b'_{#1,#2}(#3)}
\newcommand{\bs}{x}
\newcommand{\bt}{y}
\newcommand{\bx}{x}
\newcommand{\by}{y}

\newcommand{\card}{\mathtt{\#}}
\newcommand{\charpol}[1]{P_{#1}(x)}
\newcommand{\cl}[1]{[#1]}
\newcommand{\comp}[2]{[#2]_{#1}}
\newcommand{\concat}{\mathbin{{}^{\scriptscriptstyle\frown}}}

\newcommand{\D}{\Delta}
\newcommand{\DD}[2]{\D_{#1}^{#2}}
\newcommand{\dd}{d}

\newcommand{\Div}{\mathrm{Div}}
\newcommand{\dive}{\preccurlyeq}
\newcommand{\DivL}{D_{\!L}}
\newcommand{\DivR}{D_{\!R}}
\newcommand{\dual}[1]{{}^*#1}
\newcommand{\pz}{z}

\newcommand{\etc}{{\it etc.}}
\newcommand{\ev}{\rho}

\newcommand{\ff}{f}
\newcommand{\flip}[1]{\phi_{#1}}

\let\ge=\geqslant
\renewcommand{\gg}{g}

\newcommand{\ie}{{\it i.e.}}
\newcommand{\II}{I}
\let\ince=\subseteq
\newcommand{\ind}{i}

\newcommand{\Int}[1]{[\![1, #1]\!]}
\newcommand{\inv}{^{-1}}

\newcommand{\ip}{p}
\newcommand{\iq}{q}
\newcommand{\ir}{k}
\newcommand{\is}{\ell}

\newcommand{\JJ}{J}

\newcommand{\kk}{k}
\newcommand{\KK}{K}

\renewcommand{\l}{\lambda}
\let\le=\leqslant
\renewcommand{\ll}{\ell}

\newcommand{\m}{\mu}
\newcommand{\Mat}[1]{M_{#1}}
\newcommand{\Matt}[1]{M'_{#1}}
\newcommand{\Mattt}[1]{\overline M_{#1}}

\newcommand{\nbpart}[1]{p(#1)}
\newcommand{\nn}{n}

\newcommand{\nno}{{\nn-1}}

\renewcommand{\part}[2]{\{#2\}_{#1}}
\newcommand{\perm}[1]{\pi(#1)}

\newcommand{\Pw}{\mathfrak{P}}

\newcommand{\QQ}{\mathbb{Q}}

\newcommand{\resp}{{\it resp.{¨}}}
\newcommand{\rr}{r}
\newcommand{\RR}{\mathbb{R}}
\newcommand{\rrr}{r}

\let\s=\sigma
\newcommand{\set}[1]{\widetilde{#1}}
\newcommand{\sh}{\mathrm{sh}}
\renewcommand{\sp}[1]{\tau_{#1}}
\renewcommand{\ss}[1]{\sigma_{#1}}
\renewcommand{\SS}{S}
\newcommand{\sx}{x}
\newcommand{\Sym}[1]{\mathfrak{S}_{#1}}

\newcommand{\var}{u}

\newcommand{\ww}{w}

\newcommand{\xx}{x}

\newcommand{\yy}{y}

\begin{document}

\author{Patrick DEHORNOY}
\address{Laboratoire de Math\'ematiques Nicolas
Oresme UMR 6139\\ Universit\'e de Caen,
14032~Caen, France}
\email{dehornoy@math.unicaen.fr}
\urladdr{//www.math.unicaen.fr/\textasciitilde dehornoy}

\title{Combinatorics of normal sequences of
braids}

\keywords{braid group; normal form; fundamental
braid; counting; generating function}

\subjclass{20F36, 05A05}

\begin{abstract}
Many natural counting problems arise in
connection with the normal form of braids---and
seem to have not been much considered so far. Here
we solve some of them by analysing the normality
condition in terms of the associated
permutations, their descents and the
corresponding partitions. A number of different
induction schemes appear in that framework.
\end{abstract}

\maketitle

Ubiquitous and connected with a number of domains,
Artin's braid groups have received much attention
in the recent years. However, not so many works
are devoted to a purely combinatorial study of
braids, presumably because counting arguments did
not prove so far to be much helpful for
investigating braids. Nevertheless, although
braid groups are infinite, they admit several
filtrations leading to finite sets and,
therefore, to natural enumeration problems.

For each presentation of braid groups, (at least) two
natural counting poblems arise, namely, on the one hand,
counting how many braids admit an expression of
 a given length, in particular evaluating the associated
growth rate---we shall refer to this as Question¨1 in the
sequel---and, on the other hand, counting, for a given
braid, how many words represent that specific braid, a
relevant question when the number is finite,
typically when we discard the inverses of the
generators and only consider positive expressions,
\ie, when we restrict to some submonoid of~$B_\nn$---we
shall refer to that as Question¨2. 

In the case of the Artin generators¨$\ss i$, both types of
questions have been addressed, and at least partially
solved: Question¨1, actually not for¨$B_\nn$ but for the
submonoid¨$B_\nn^+$ of¨$B_\nn$ generated by the¨$\ss i$'s, was
investigated in~\cite{Xua}, and completely solved
in~\cite{Bro}. As for Question~2, it is natural in this context to
address it for the particular elements¨$\D_n^d$, where
$\D_n$ is Garside's fundamental braid \cite{Gar}. It was
investigated and solved for~$n = 3$ in~\cite{CrH}.

In this paper, we address Question~1 for another natural 
generating set, namely the so-called simple braids, also called
the Garside generators below~\cite{Gar}. 
These generators, which are the divisors  of~$\D_n$ in 
the monoid~$B_\nn^+$, are in one-to-one
correspondence with permutations of $n$~objects,
and they give rise to a remarkable unique
decomposition for each braid, usually called its normal
form \cite{Dlg, Adj, Eps, ElM}. Because of its
uniqueness and of its many nice properties,
expressed in particular in the existence of a bi-automatic 
structure, the normal form
of braids is the preferred way of specifying
braids in many recent developments, in particular
those of algorithmic or cryptographical nature \cite{Geb, KoL}.

In the case of the Garside generators and the associated
normal form(s), what we called Question¨2 above is not
relevant, as each braid has one unique distinguished
decomposition in terms of the generators, and counting the non-normal expressions would appear artificial. But Question¨1, namely counting the
number of braids with a normal form of a given length, is
quite natural. The question was briefly considered
by R.\,Charney in¨\cite{Chb}: she observed that, because
normal words can be recognized by a finite state automaton,
the number of braids with length¨$\dd$ obeys a linear
induction rule, and the associated generating function is
 rational. Explicit values are given in the case of
$3$¨and $4$¨strand braids.

The aim of this paper is to go further in the investigation
of counting problems connected with the normal form of
braids. We consider the case of the monoid¨$B_\nn^+$, and
mainly study the number of positive $n$¨strand braids
with a normal form of length (at most)¨$\dd$, \ie, we
consider for Garside generators the problem solved
in¨\cite{Bro} in the case of Artin generators. The main
difference is that, in the case of
Garside generators, the length is no longer an additive
parameter: for instance, multiplying two simple braids may
result in a simple braid, \ie, multiplying two braids of
length¨$1$ may result in a braid of length¨$1$. This makes
the current study much more uneasy.

Let $\bb\nn\dd$ denote the number of positive $n$¨strand
braids with a normal form of length at most¨$\dd$, \ie, the
number of divisors of¨$\DD\nn\dd$ in¨$B_\nn^+$. We establish
various results about the numbers¨$\bb\nn\dd$,
and about the connected numbers¨$\bbb\nn\dd\bs$ that count,
for $\bs$ a simple braid, the positive $n$¨strand braid
with a normal form of length at most¨$\dd$ whose
$\dd$th factor is precisely¨$\bs$. Two types of results are established, 
namely results for fixed braid index¨$n$, and results for fixed degree¨$\dd$.
When $n$ is fixed and $\dd$ varies, as was recalled above, 
the numbers¨$\bb\nn\dd$ and $\bbb\nn\dd\bs$ obey a linear
induction rule associated with a certain $n! \times n!$
incidence matrix¨$\Mat\nn$. Here we show that $\Mat\nn$
can be replaced with a smaller matrix of size $p(\nn)
\times p(\nn)$, where $p(\nn)$ is the number of partitions
of¨$\nn$. The result relies on analysing the descents
of the permutations associated with simple braids, and it
is connected with a classical result by Solomon¨\cite{Sol}.
It is then easy to deduce the numerical value of¨$\bb\nn\dd$
for small¨$\nn, \dd$, as well as explicit formulas, at least
for $n \le 4$. We are also led to several
conjectures about the eigenvalues of the matrix¨$\Mat\nn$
that seem to have never been considered so far. The most
puzzling one claims that the characteristic polynomial
of¨$\Mat{\nno}$ divides that of¨$\Mat\nn$. It holds at
least for $\nn \le 10$.

When $\dd$ is fixed and $\nn$ varies, quite different
induction rules appear. Everything is trivial for $\dd =
1$, and an explicit formula for¨$\bb\nn2$ can be deduced from
the results of¨\cite{CSV, CSV2}. It seems difficult to
go further in general, but new results (and new
induction schemes) appear when we consider the numbers
$\bbb\nn\dd{\D_{\nn-\rrr}}$ with $1 \le \rr \le \nn$,
typically in the (non-trivial) case $\rr = 1$, and, more
generally, when $\rr$ is fixed. In particular, we
obtain explicit values for $\bbb\nn3{\D_{\nno}}$,
$\bbb\nn3{\D_{\nn-2}}$, and $\bbb\nn4{\D_{\nno}}$.

The specific questions investigated in this paper, in
particular that of the value
of¨$\bbb\nn\dd{\D_{\nn-\rrr}}$, arose in¨\cite{Dhh}. There
exists a distinguished linear ordering of braids¨$<$ relying
on the notion of a
$\s$-positive braid word¨\cite{Dgr}, and the aim
of¨\cite{Dhh} is to develop a new approach to that ordering
based on the study of its connection with the Garside
structure. It turns out that certain parameters describing the
restriction of¨$<$ to positive $n$-braids of degree at
most¨$\dd$ can be expressed in terms of the
numbers¨$\bbb\nn\dd{\D_\rr}$, an initial
motivation for our current study of these numbers. However,
we think that the formulas and methods developed in the
current paper go beyond the above specific applications. In
particular, the great diversity of the induction
schemes appearing in connection with various
specializations of the general problem is remarkable. At the
least, the current study should demonstrate the richness of
the combinatorics underlying the normal form of braids.

Still other presentations of the braid groups are known, in
particular the one involving the so-called dual
monoid¨\cite{Bes, BKL}, which gives rise to an alternative
Garside structure, and, therefore, to an alternative
normal form analogous to that considered here, where the
role of simple braids is played by elements that are in
one-to-one correspondence with non-crossing partitions.
All questions considered in the current paper could be
similarly addressed for the dual structure, and, more
generally, for the many presentations of¨$B_\nn$ known to
date. Similarly, Artin's braid groups~$B_\nn$ belong to larger
families of groups, typically Artin-Tits groups of spherical
type and, more generally, Garside groups~\cite{Dfx, Dgk, McC}. 
Once again, all
questions considered here extend to such frameworks
naturally. However, mainly because of the specific
applications mentioned above, we find it interesting to
consider here the specific framework of braids and
permutations, and we leave the extensions for further
investigation.

The paper is organized as follows. Section¨\ref{S:Back}
sets the framework and the basic definitions. In
Section¨\ref{S:Adj} we introduce the incidence
matrix¨$\Mat\nn$ that controls the sequences¨$\bb\nn\dd$ for
fixed¨$\nn$ and show how to reduce their size
from¨$\nn!$ to¨$2^\nno$. In Section¨\ref{S:Part}, we
show how to further reduce the size to¨$p(\nn)$, and solve
the induction for small values of¨$\nn$. Finally, in
Section¨\ref{S:Deg}, we turn to the cases when the
degree is fixed and the braid index varies. 

\section*{Acknowledgment}

The author thanks C.\,Hohlweg, F.\,Hivert and
J.C.\,Novelli for interesting discussions about the topics
investigated in this paper, in particular Remark¨\ref{R:Hohlweg}.

  \section{Background and preliminary results}
\label{S:Back}

Our notation is standard, and we refer to textbooks
like¨\cite{Bir} or¨\cite{Eps} for basic results about
braid groups. We recall that the $\nn$¨strand braid group¨$B_\nn$
is defined for $\nn \ge 1$ by the presentation
\begin{equation} \label{E:Present}
B_\nn = \left<\ss1, \dots, \ss{\nno} \,;\,
\begin{array}{cl}
\ss i \ss j = \ss j \ss i 
&\text{\quad for $\vert
i - j\vert \ge 2$}\\
\   \ss i \ss j \ss i  = \ss j \ss i \ss j
&\text{\quad for $\vert i - j\vert = 1$}
\end{array} \right>.
\end{equation}
So, $B_1$ is a trivial group¨$\{1\}$, while $B_2$
is the free group generated by¨$\ss1$. The elements of¨$B_\nn$
are called $\nn$¨strand braids, or simply {\it $\nn$-braids}.
We use $B_\infty$ for the group generated by an infinite
sequence of¨$\ss i$'s subject to the relations
of¨\eqref{E:Present}, \ie, the direct limit of all¨$B_\nn$'s
under the inclusion of¨$B_\nn$ into¨$B_{n+1}$.

By definition, every $\nn$-braid¨$\bx$ admits (infinitely
many) expressions in terms of the generators¨$\ss i$, $1
\le i < n$. Such a expression is called an $\nn$¨strand {\it
braid word}. Two braid words¨$\ww, \ww'$ representing the
same braid are said to be {\it equivalent}; the braid represented 
by a braid word¨$\ww$ is denoted¨$\cl\ww$. 

It is standard to associate with every $\nn$¨strand braid
word¨$\ww$ an $\nn$¨strand braid {\it diagram} by stacking 
elementary diagrams associated with the successive letters 
according to the rules
\begin{center}
\begin{picture}(83,13)(0,0)
\put(10,0){\includegraphics{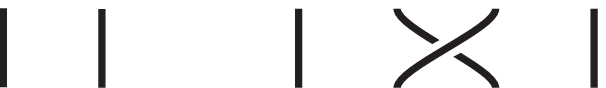}}
\put(-4,3){$\ss i\mapsto$}
\put(9.5,9){$1$}
\put(19.5,9){$2$}
\put(50,9){$i$}
\put(58,9){$i{+}1$}
\put(27,3){\dots}
\put(77,3){\dots}
\end{picture}
\end{center}
\begin{center}
\begin{picture}(83,10)(0,0)
\put(10,0){\includegraphics{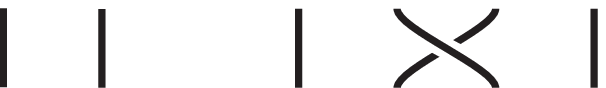}}
\put(-4,3){$\ss i\inv \mapsto$}
\put(27,3){\dots}
\put(77,3){\dots}
\end{picture}
\end{center}
Then two braid words are
equivalent if and only if the diagrams they encode are the
projections of ambient isotopic figures in¨$\RR^3$, \ie, one
can deform one diagram into the other without allowing the
strands to cross or moving the endpoints. 

\subsection{The monoid¨$B_\nn^+$ and the braids¨$\D_\nn$}

Let $B_\nn^+$ be the monoid admitting
the presentation¨\eqref{E:Present}. The elements of¨$B_\nn^+$
are called {\it positive} $\nn$-braids.

\begin{defi}
For $\bx, \by$ in¨$B_\nn^+$, we say that $\bx$ is a {\it
left divisor} of¨$\by$, denoted $\bx \dive \by$, or,
equivalently, that $\by$ is a {\it right multiple}
of¨$\bx$, if $\by = \bx \pz$ holds for some¨$\pz$
in¨$B_\nn^+$. We denote by¨$\Div(\by)$ the (finite) set of
all left divisors of¨$\by$ in¨$B_\nn^+$.
\end{defi}

As $B_\nn^+$ is not commutative for $\nn \ge 3$, there are
the symmetric notions of a right divisor and a left
multiple---but we shall mostly use left divisors here.
Note that $\bx$ is a (left) divisor of¨$\by$ in the sense
of¨$B_\nn^+$ if and only if it is a (left) divisor in the
sense of¨$B_\infty^+$, so there is no need to
always specify the index¨$\nn$.

Wih respect to left divisibility, $B_\nn^+$ has the structure
of a lattice \cite{Gar}: any two positive $\nn$-braids¨$\bx,
\by$ admit a greatest common left divisor, denoted
$\gcd(\bx, \by)$, and a least common right multiple. A
special role is played by the lcm of the elements¨$\ss1$,
\dots, $\ss{\nno}$, traditionally denoted¨$\D_\nn$, which we
recall is inductively defined by
\begin{equation} \label{E:Delta}
\D_1 = 1, \qquad
\D_\nn = \ss1 \ss2 \dots \ss{\nno} \, \D_{\nno}.
\end{equation}
It is well known that $\D_\nn^2$ belongs to the
centre of¨$B_\nn$ (and even generates it for $\nn \ge 3$), and
that the inner automorphism¨$\flip n$ of¨$B_\nn$ corresponding
to conjugation by¨$\D_\nn$ exchanges¨$\ss i$ and¨$\ss{n-i}$
for $1 \le i \le \nno$. 

\subsection{The normal form}

In¨$B_\nn^+$, the left and the right divisors of¨$\D_\nn$
coincide, and they make a finite sublattice of¨$(B_\nn^+,
\dive)$ with $\nn!$¨elements. These braids will be called
{\it simple} in the sequel. Geometrically, simple braids
are those positive braids that can be represented by a
braid diagram in which any two strands cross at most once.

For each positive $\nn$-braid¨$\bx$ distinct of¨$1$,
the simple braid $\gcd(\bx, \D_\nn)$ is the maximal
simple left divisor of¨$\bx$, and we obtain a
distinguished expression $\bx = \sx_1 \bx'$
with¨$\sx_1$ simple. By decomposing¨$\bx'$ in the
same way and iterating,  we obtains the so-called
normal expression¨\cite{ElM, Eps}. 

\begin{defi}
A sequence $(\sx_1, \dots, \sx_\dd)$ of simple
$\nn$-braids is said to be {\it normal} if, for
each¨$\kk$, one has $\sx_\kk =
\gcd(\D_\nn, \sx_\kk \dots \sx_\dd)$.
\end{defi}

Clearly, each positive braid admits a unique normal
expression. It will be convenient here to consider the
normal expression as unbounded on the right by
completing it with as many trivial factors¨$1$ we need. In
this way, we can speak of the {\it $\dd$th factor} (in the
normal form) of¨$\bx$ for each positive braid¨$\bx$. We
say that a positive braid has {\it degree¨$\dd$} if $\dd$ is
the largest integer such that the $\dd$th factor of¨$\bx$ is
not¨$1$. It is well known that the positive $\nn$-braids of
degree at most¨$\dd$ coincide with the (left or oright)
divisors of¨$\DD\nn\dd$. 

The only properties of the normal form we shall use here
are as follows:

\begin{lemm} \cite{ElM} \label{L:Normal}
Assume that $(\sx_1, \dots, \sx_\dd)$ is a sequence
of simple
$\nn$-braids. Then the following are equivalent:

$(i)$ The sequence $(\sx_1, \dots, \sx_r)$ is normal;

$(ii)$ For $1 \le \kk < \dd$, the sequence
$(\sx_\kk, \sx_{\kk+1})$ is normal;

$(iii)$ For $1 \le \kk < \dd$, every¨$\ss i$
dividing¨$\sx_{\kk+1}$ on the left divides¨$\sx_\kk$
on the right.
\end{lemm}

\begin{defi}
For¨$\sx$ a simple $\nn$-braid, we define
$\DivL(\sx)$ (\resp $\DivR(\sx)$) to be the set of
all¨$i$'s such that $\ss i$ is a left (\resp
right) divisor of¨$\sx$.
\end{defi}

The example of¨$\ss2$ and $\ss2\ss1\ss3\ss2$,
for which both¨$\DivL$ and¨$\DivR$ is¨$\{2\}$,
shows that these sets do not determine a
simple braid. However, as far as normal sequences are
concerned, they contain all needed
information, as Lemma¨\ref{L:Normal}
can be restated as:

\begin{lemm} \label{L:NormalBis}
A sequence of simple $\nn$-braids $(\sx_1, \dots,
\sx_\dd)$ is  normal if and only if, for each¨$\kk <
\dd$, we have $\DivR(\sx_\kk)
\supseteq \DivL(\sx_{\kk+1})$.
\end{lemm}

\subsection{Connection with permutations}

Everywhere in the sequel, we write¨$\Int\nn$ for $\{1,
\dots, \nn\}$. By mapping¨$\ss i$ to the transposition¨$(i,
i+1)$, one defines a surjective homomorphism¨$\pi$
of¨$B_\nn$ onto the symmetric group¨$\Sym\nn$. The
restriction of¨$\pi$ to simple braids is a bijection:
for every permutation¨$f$ of¨$\Int\nn$, there exists exactly
one simple braid¨$\bs$ satisfying $\perm\bs = f$.

The Exchange Lemma for Coxeter groups connects
the sets¨$\DivL(\sx)$ and¨$\DivR(\sx)$ with the
permutation associated with¨$\sx$ and their descents. 
For¨$\ff$ a permutation, use $\ell(\ff)$ for the minimal 
number of factors occurring in a decomposition of¨$\ff$ 
as a product of transpositions. The precise statement is

\begin{lemm} \label{L:Exchange}
Let $\bs$ be a simple $\nn$-braid¨$\bs$. For $1 \le i
< \nn$, the following are equivalent:

$(i)$ The braid¨$\ss i$ is a left divisor of¨$\bs$
in¨$B_\nn^+$, \ie, $i$ belongs to¨$\DivL(\bs)$; 

$(ii)$ The strands starting at positions¨$i$ and¨$i+1$
cross in any (positive) diagram for¨$\bs$;

$(iii)$ We have $\pi(\bs)\inv(i) > \pi(\bs)\inv(i+1)$;

$(iv)$ We have $\ell(\pi(\ss i\bs)) < \ell(\perm\bs)$, \ie, $i$ is
a descent of¨$\perm\bs\inv$.\\
Symmetrically, the following are equivalent:

$(i')$ The braid¨$\ss i$ is a right divisor of¨$\bs$
in¨$B_\nn^+$, \ie, $i$ belongs to¨$\DivR(\bs)$;  

$(ii')$ The strands finishing at positions¨$i$ and¨$i+1$
cross in any (positive) diagram for¨$\bs$;

$(iii')$ We have $\pi(\bs)(i) > \pi(\bs)(i+1)$;

$(iv')$ We have $\ell(\pi(\bs\ss i)) < \ell(\perm\bs)$, \ie, $i$ is
a descent of¨$\perm\bs$.
\end{lemm}

So, for¨$\bs$ a simple braid, the indices¨$i$ such that $\ss i$ is a right
divisor of¨$\bs$ are the descents of the associated permutation¨$\perm\bs$,
while those such that $\ss i$ is a left divisor of¨$\bs$ are the
descents of¨$\perm\bs\inv$.

\subsection{The numbers $\bb\nn\dd$ and $\bbb\nn\dd\bs$}

Our aim in this paper is to solve various counting problems
involving the normal form of positive braids. The main
numbers we investigate are as follows:

\begin{defi}
For $\nn, \dd \ge 1$, we denote by¨$\bb\nn\dd$ the number
of positive $\nn$¨strand braids of degree at most¨$\dd$,
\ie, the number of divisors of¨$\DD\nn\dd$ in the braid
monoid.
\end{defi}

By Lemma¨\ref{L:NormalBis}, $\bb\nn\dd$ is the number of
normal sequences of length¨$\dd$, \ie, the number of
sequences $(\sx_1, \dots, \sx_\dd)$ where all¨$\sx_\kk$ are
simple braids and $\DivL(\sx_\kk) \supseteq
\DivR(\sx_{\kk+1})$ holds for $\kk < \dd$. By
Lemma¨\ref{L:Exchange}, it is also the number of sequences
of permutations $(\ff_1, \dots, \ff_\dd)$ such that, for
each¨$\kk < \dd$, the descents of¨$\ff_{\kk+1}\inv$ are
included in those of¨$\ff_\kk$.

For¨$\dd = 1$, the bijection between simple $\nn$¨strand
braids and permutations of¨$\Int\nn$ immediately gives
\begin{equation}
\bb\nn1 = n!,
\end{equation}
which implies for all¨$\nn, \dd$
\begin{equation} \label{E:Upper}
\bb\nn\dd \le (\nn!)^\dd.
\end{equation}

In the sequel, we shall have to count normal sequences
satisfying some constraints. So we introduce one more
notation.

\begin{defi}
For $\nn, \dd \ge 1$ and $\bs$ a simple $\nn$-braid, we
denote by¨$\bbb\nn\dd\bs$ the number of positive
$\nn$¨strand braids of degree at most¨$\dd$ with $\dd$th
factor equal to¨$\bs$.
\end{defi}

In other words, $\bbb\nn\dd\bs$ is the number of normal
sequences of the form $(\bs_1, \dots, \bs_{\dd-1}, \bs)$.
Some connections are obvious:

\begin{prop} \label{P:Total}
For all¨$\nn, \dd$, we have
\begin{equation} \label{E:Last}
\bb\nn\dd = \sum_{\mbox{\Small\rm $\bs$ simple}} \bbb\nn\dd\bs =
\bbb\nn{\dd+1}1.
\end{equation}
\end{prop}

\begin{proof}
The first equality is obvious. The second one follows from
the fact that $(\bs_1, \dots, \bs_\dd)$ is normal if and
only if $(\bs_1, \dots, \bs_\dd, 1)$ is: indeed, $1$ has no
left divisor but itself, so, by Lemma¨\ref{L:Normal}, every
sequence $(\bs, 1)$ is normal.
\end{proof}

\section{Adjacence matrices}
\label{S:Adj}

In this section, we study the numbers $\bb\nn\dd$ and
$\bbb\nn\dd\bs$ when $\nn$ is fixed and $\dd$ varies. By
Lemma¨\ref{L:Normal}, normal sequence of simple braids are
characterized by a purely local criterion that only involves
adjacent entries. It follows that the set of all normal
sequences can be recognized by a finite state automaton
\cite{Eps}, and, as a consequence, the associated counting
numbers obey a linear induction rule specified by a certain
incidence matrix¨\cite{Eps2}. In this section, we define
the matrix involved in the current situation, and show how
its size, which is originally¨$\nn!$, can be lowered
to¨$2^\nno$.

\subsection{Enumeration of simple braids}

Below we consider matrices whose entries are indexed
by simple braids (or, equivalently, permutations). Fixing an
enumeration of simple braids is not important at a
conceptual level, but this is necessary when the objects
are to be specified explicitly. We shall use the
restriction of the canonical linear ordering of braids
denoted¨$<^{\scriptscriptstyle\phi}$
in¨\cite{Dgr}---which gives for each¨$\nn$ a well-ordering of
ordinal type¨$\omega^{\omega^{\nn-2}}$ on¨$B_\nn^+$. 
The corresponding increasing enumeration of simple
$\nn$-braids can be constructed directly using induction
on¨$\nn$. We start from the following easy remark:

\begin{lemm} \label{L:Factor}
For $1 \le i \le \nn$, write $\ss{i, \nn}$ for $\ss i
\ss{i+1} \dots \ss{\nno}$ (so that $\ss{\nn, \nn}$
is¨$1$). Then every simple $\nn$-braid¨$\bs$ admits a unique
decomposition $\bs = \ss{i, \nn} \, \bt$ with $1
\le i \le \nn$ and $\bt$ a simple $(\nno)$-braid.
\end{lemm}
 
\begin{figure} [htb]
\begin{picture}(35,20)(0, 0)
\put(0,0){\includegraphics{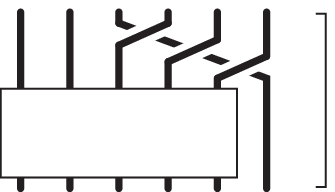}}
\put(34,9){$\bs$}
\put(11.5,5){$\bt$}
\put(12,19.5){$i$}
\end{picture}
\caption{\smaller Proof of Lemma¨\ref{L:Factor}}
\label{F:Factor}
\end{figure}

\begin{proof}
(Figure¨\ref{F:Factor})
Let $i = \perm\bs(\nn)$. Then we can realize¨$\bs$ by a
diagram in which the $i$th strand is first sent to the
rightmost position, and it remains a simple $(\nno)$-braid.
Conversely, we have $i = \perm{\ss{i, \nn}\bt}(\nn)$, so the
decomposition is unique.
\end{proof}

\begin{defi}
We inductively define an enumeration¨$S_\nn$ of simple
$\nn$-braids by
\begin{equation}
S_1:= (1), \quad
S_\nn:= S_{\nno} \concat \ss{\nno, \nn} S_{\nno} \concat
\dots \concat \ss{1, \nn} S_{\nno},
\end{equation}
where $\concat$ stands for list concatenation, and $\bs S$ is
the list obtained from¨$S$ by multipliying all entries
by¨$\bs$ on the left. The $\kk$th element in¨$\bigcup_\nn
S_\nn$ is denoted¨$\sp\kk$.
\end{defi}

The first¨$\sp\kk$'s are, in increasing order,
$$\sp1 = 1\ , 
\sp2 = \ss1\ , 
\sp3 = \ss2\ , 
\sp4 = \ss2\ss1\ ,
\sp5 =  \ss1\ss2\ , 
\sp6 = \ss1\ss2\ss1\ , 
\sp7 = \ss3, \dots$$
Lemma¨\ref{L:Factor} guarantees that all simple braids
occur in the above enumeration. Note that, for every¨$\nn$,
we have $\D_\nn = \sp{\nn!}$.

It is easy to check that
the ordering of simple braids we use corresponds to a
reversed antilexicographic ordering of the inverses of the
associated permutations: $\bs$ occurs before¨$\bt$ if and
only if we have $\perm\bs\inv < \perm\bt\inv$, where 
$\ff < \gg$ is said to hold if we have $\ff(i) > \gg(i)$ for
the largest¨$i$ for which $\ff$ and $\gg$ do not agree.

\subsection{The matrix¨$\Mat\nn$}

Everywhere in the sequel, we write $(M)_{\xx, \yy}$ for the
$(\xx, \yy)$-entry of a matrix¨$M$.

\begin{defi} 
For $\nn \ge 1$, we define $\Mat\nn$ to be the $\nn! \times
\nn!$ matrix satisfying
$$
(\Mat\nn)_{\kk, \ll} = 
\begin{cases}
1 & \mbox{if $(\sp\kk, \sp\ll)$ is normal,}\\
0 & \mbox{otherwise.}
\end{cases}$$
\end{defi}

Instead of referring to integer entries, it will be often
convenient to think of the entries of¨$\Mat\nn$ as directly
indexed by simple braids; for $\bs, \bt$ simple braids, we
simply write $(\Mat\nn)_{\bs, \bt}$ for the
corresponding entry.

\begin{exam}
The first 3 matrices¨$\Mat\nn$ are
$$
\Mat1 = 1, \quad
\Mat2 = 
\left(\begin{matrix}1&0\\1&1\end{matrix}\right), \quad
\Mat3 = \left(
\begin{matrix}
1&0&0&0&0&0\\
1&1&0&0&1&0\\
1&0&1&1&0&0\\
1&1&0&0&1&0\\
1&0&1&1&0&0\\
1&1&1&1&1&1
\end{matrix}
\right).$$
\end{exam}

The construction of the matrix¨$\Mat\nn$ immediately
implies the following results:

\begin{lemm} \label{L:Stack}
$(i)$ The first column and the last row of¨$\Mat\nn$
contain only¨$1$'s; the first row, its first entry
excepted, and the last column, its last entry excepted, 
contain only¨$0$'s.

$(ii)$ The first $(\nno)!$ columns of¨$\Mat\nn$ consist of
$\nn$¨stacked copies of¨$\Mat{\nno}$. 

$(iii)$ If $\DivR(\sp\kk) = \DivR(\sp{\kk'})$ holds, then
the $\kk$th and the $\kk'$th rows in¨$\Mat\nn$ coincide.
Similarly, if $\DivL(\sp\ll) = \DivL(\sp{\ll'})$ holds, then
the $\ll$th and the $\ll'$th columns in¨$\Mat\nn$ coincide.
\end{lemm}

\begin{proof}
$(i)$ By  construction, we have $\sp1 = 1$ and
$\sp{\nn!} = \DD\nn{}$. Now $(\bs, 1)$ is always normal,
and so is $(\DD\nn{}, \bs)$. On the other hand, $(1 , \bs)$
is normal only for $\bs = 1$, and $(\bs, \DD\nn{})$ is
normal only for $\bs = \DD\nn{}$.

$(ii)$ Assume $\kk \le (\nno)!$ and $\kk' = \kk + (\nn-i)
\cdot (\nno)!$. Our enumeration of simple braids implies
$\sp{\kk'} = \ss{i, \nn} \sp\kk$. Then
Figure¨\ref{F:Factor} makes the equality
$\DivR(\sp{\kk'}) \cap \Int{\nno} = \DivR(\sp\kk)$ clear.
For every simple $(\nno)$-braid¨$\bt$, the
set¨$\DivL(\bt)$ is included in¨$\Int{\nno}$, and it
follows that $(\sp{\kk'}, \bt)$ is normal if and only
if¨$(\sp\kk, \bt)$ is. In other words, we have 
$(\Mat\nn)_{\kk', \ll} = (\Mat\nn)_{\kk, \ll}$ for $\ll \le
(\nno)!$.

$(iii)$ By Lemma¨\ref{L:Normal}, the value of
$(\Mat\nn)_{\kk, \ll}$ only depends on¨$\DivR(\sp\kk)$ and
on¨$\DivL(\sp\ll)$.
\end{proof}

The connection between the numbers¨$\bbb\nn\dd\bs$ and the
matrix¨$\Mat\nn$ is straightforward:

\begin{lemm} \label{L:Comput}
For every simple¨$\bt$ and every
$\dd \ge 1$, we have
\begin{equation} \label{E:Number}
\bbb\nn\dd{\bt} = ((1, 1, \dots, 1) \,
\Mat\nn\!\!\!{}^{\dd-1})_\bt.
\end{equation}
\end{lemm}

\begin{proof}
Induction on¨$\dd$. For $\dd=1$, and for each simple
$\nn$-braid¨$\bs$, there is exactly one braid of
degree at most¨$1$ whose first factor is¨$\bt$,
namely $\bt$¨itself, and we have $\bbb\nn1\bt = 1$.
Assume $\dd \ge 2$. By Lemma¨\ref{L:Normal}, $(\bs_1,
\dots, \bs_{\dd-1}, \bt)$ is normal if and only if $(\bs_1,
\dots, \bs_{\dd-1})$ and $(\bs_{\dd-1}, \bt)$ are normal,
so we get 
\begin{equation*} \label{E:Comput}
\bbb\nn\dd{\bt} 
= \sum_{(\bs, \bt)\mbox{ \Small normal}}
\bbb\nn{\dd-1}{\bs} = \sum_\bs
\bbb\nn{\dd-1}{\bs}
\, (\Mat\nn)_{\bs, \bt},
\end{equation*}
and \eqref{E:Number} follows inductively. 
\end{proof}

\begin{rema}
As last row of¨$\Mat\nn$ is $(1, \dots, 1)$, we have
$(1, \dots, 1) = (0, \dots, 0, 1) \Mat\nn$, and we can
replace¨\eqref{E:Number} with
\begin{equation} \label{E:NumberBis}
\bbb\nn\dd{\bt} = ((0, \dots, 0, 1) \,
\Mat\nn\!\!\!{}^{\dd})_\bt.
\end{equation}
\end{rema}

\begin{exam}
Using the value of¨$\Mat2$, we immediately
find $\bbb2\dd1 = \dd$, $\bbb2\dd{\ss1} = 1$, as
could be expected: there are $\dd+1$ braids of degree
at most¨$\dd$, namely the braids¨$\ss1^\kk$ with
$\kk < \dd$, whose $\dd$th factor is¨$1$, and
$\ss1^\dd$, whose $\dd$th factor is¨$\D_2$, \ie,¨$\ss1$. 
\end{exam}

The computation for $\nn \ge 3$ is more complicated, and we
postpone it. For the moment, we just point that, as the
numbers¨$\bbb\nn\dd\bs$ obey the linear
recurrence¨\eqref{E:Number}, standard arguments imply  that
they can be expressed in terms of the eigenvalues
of¨$\Mat\nn$:

\begin{prop} \label{P:Expansion}
Let $\ev_1$, \dots, $\ev_\ir$ be the non-zero
eigenvalues of¨$\Mat\nn$. Then,  for
each simple $\nn$-braid¨$\bs$, there
exist polynomials $P_1, \dots, P_\ir$ with $\deg(P_\ind)$
at most the multiplicity of¨$\ev_\ind$ for¨$\Mat\nn$ such
that, for each¨$\dd \ge 0$, we have
\begin{equation} \label{E:General}
\bbb\nn\dd\bs = P_1(\dd) \ev_1^\dd + \cdots +
P_\ir(\dd) \ev_\ir^\dd.
\end{equation} 
\end{prop}

\begin{coro}
For all¨$\nn, \bs$, the generating function of the
numbers¨$\bbb\nn\dd\bs$'s with respect to¨$\dd$ is rational.
\end{coro}

\subsection{Reducing the size}

The size¨$\nn!$ of the incidence¨$\Mat\nn$ is uselessly
large, and we shall see now how to lower it. This will be
done in two steps. The first one relies on the fact,
pointed out in Lemma¨\ref{L:Stack}$(iii)$, that many columns
in¨$\Mat\nn$ are equal. For subsequent use, it will be
useful to introduce a new sequence of numbers:

\begin{defi}
For $\II, \JJ \ince \Int{\nno}$, we denote
by $\aa\nn\II\JJ$ (\resp $\aas\nn\II\JJ$) the
number of simple $\nn$-braids satisfying $\DivL(\bs) =
\II$ (\resp $\DivL(\bs) \supseteq \II$) and
$\DivR(\bs) \supseteq \JJ$.
\end{defi}

\begin{lemm} \label{L:Reduction}
For $\nn \ge 1$, let $\Matt\nn$ be the $2^\nno \times
2^\nno$ matrix with entries indexed by subsets
of¨$\Int{\nno}$ defined by $(\Matt\nn)_{\II, \JJ} = 
\aa\nn\II\JJ$. Then the characteristic polynomials
of¨$\Matt\nn$ and $\Mat\nn$ coincide up to a power of¨$\xx$,
and, for every simple¨$\bt$ with $\DivL(\bt) = \JJ$ and
every $\dd \ge 1$, we have
\begin{equation} \label{E:Numberr}
\bbb\nn\dd{\bt} = ((1, 1, \dots, 1) \,
\Matt\nn\!{}^{\dd-1})_\JJ.
\end{equation}
\end{lemm}

\begin{proof}
Gathering the columns corresponding to simples with
the same¨$\DivL$¨set and summing the corresponding
lines amounts to replacing¨$\Mat\nn$ with a similar
matrix of the form
$\left(\begin{matrix}
\Matt\nn&0\\\dots&0 \end{matrix}\right)$, so the
result about the characteristic polynomial is clear.

As for the value of¨$\bbb\nn\dd\bt$, the argument
is similar to that for Lemma¨\ref{L:Comput}. The
induction starts as $\aa\nn\II{\Int{\nno}} =
1$ holds for each¨$\II$. For the general step, we
find
\begin{align*} \label{E:Comput}
\bbb\nn\dd{\bt} 
&= \sum_{(\bs, \bt)\mbox{ \Small normal}}
\bbb\nn{\dd-1}{\bs} 
= \sum_{\DivR(\bs) \supseteq\JJ}
\bbb\nn{\dd-1}{\bs} 
= \sum_\II \sum_{
\underset{\scriptstyle \DivL(\bs) = \II}
{\scriptstyle \DivR(\bs) \supseteq\JJ}}
\bbb\nn{\dd-1}{\bs} \\
&= \sum_\II  \bbbb\nn{\dd-1}{\II} \aa\nn\II\JJ
= \sum_\II  ((1, \dots, 1)
\Matt\nn{}^{\dd-2})_{\II} (\Matt\nn)_{\II, \JJ}
= ((1, \dots, 1)
\Matt\nn{}^{\dd-1})_{\JJ},
\end{align*}
where $\bbbb\nn\dd\II$ denotes the common
value of $\bbb\nn\dd\bs$ for $\bs$ with¨$\DivL(\bs) =
\II$. 
\end{proof}

For $\nn = 3$, and using the enumeration $\emptyset$,
$\{1\}$, $\{2\}$, $\{1,2\}$ that is induced by our
enumeration of simple braids, we obtain 
$\Matt3 = \left(\begin{matrix}
1&0&0&0\\2&1&1&0\\
2&1&1&0\\
1&1&1&1
\end{matrix}\right)$. Observe that the second
and third columns in¨$\Matt3$ coincide, which suggests a
further reduction step.

\section{Partitions associated with a simple braid}
\label{S:Part}

We can indeed reduce the size of the matrices once more: 
we can replace the
incidence matrix¨$\Matt\nn$ with a new matrix¨$\Mattt\nn$,
whose size is¨$p(\nn)$, the number of partitions of¨$\nn$.
Here the result is deduced from elementary remarks
about simple braids (or, equivalently, about
permutations); it can also be deduced from classical
results about Solomon's algebra of descents---and 
therefore extends to all Artin--Tits groups of
spherical type. 

\subsection{Computation of¨$\aas\nn\II\JJ$}

We shall start from an explicit determination of the value
of the numbers¨$\aas\nn\II\JJ$ in terms of the block
compositions of¨$\II$ and¨$\JJ$. We first recall the notions
of composition and partition. 

\begin{defi}
For $\II \ince \Int{\nno}$, the {\it
$\nn$-composition}¨$\comp\nn\II$ of~$\II$ is defined to be
the sequence $(\ip_1, \ip_2-\ip_1,
\dots,
\ip_\ir-\ip_{\ir-1})$, where
$\ip_1, \dots, \ip_\ir$ is the increasing
enumeration of $\Int\nn \setminus \II$. The {\it
$\nn$-partition}¨$\part\nn\II$ of~$\II$ is the
non-increasing rearrangement of¨$\comp\nn\II$.
\end{defi}

\begin{exam}
The composition of~$\II$ consists of the sizes of the blocks
of adjacent elements in~$\II$, augmented by¨$1$: for
instance, the $10$-composition of
$\{1,2,4,5,6,9\}$ is¨$(3, 4, 1, 2)$, as the blocks are  of
size¨$2, 3, 0, 1$. The $10$-partition of~$\II$ is therefore
$(4, 3, 2, 1)$. Note that the $\nn$-composition of~$\II$
determines~$\II$, but its $\nn$-partition does not.
\end{exam}

The geometric observation is the following one:

\begin{lemm} \label{L:Sequence}
Assume $\JJ \ince \Int{\nno}$, and let
$(\iq_1, \dots, \iq_\is)$ be the $\nn$-composition
of¨$\JJ$. For $\bs$ a simple $\nn$-braid, define
$\ff_\bs : \Int\nn \to \Int\is$ by $\ff_\bs(i) = j$
if the $i$th strand of¨$\bs$ finishes in the
$j$th block of¨$\JJ$. Then $\bs \mapsto \ff_\bs$
establishes a bijection between the simple
$\nn$-braids¨$\bs$ that satisfy
$\DivR(\bs) \supseteq \JJ$ and the functions
from¨$\Int\nn$ to¨$\Int\is$. Moreover, for
$\DivR(\bs) \supseteq \JJ$, we have, for every¨$i$, 
\begin{equation} \label{E:Equiv}
i \in \DivL(\bs) 
\quad  \Longleftrightarrow \quad
\ff_\bs(i) \ge \ff_\bs(i+1).
\end{equation}
\end{lemm} 

\begin{figure} [htb]
\begin{picture}(92,44)(0, 0)
\put(0,4){\includegraphics{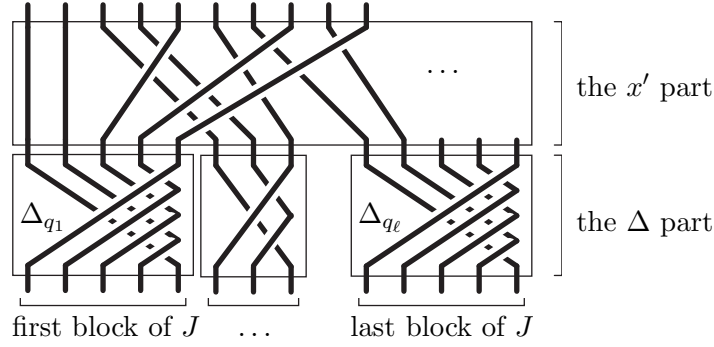}}
\put(75,32){the $\bs'$ part}
\put(75,14){the $\D$ part}
\put(1,15){$\D_{\iq_1}$}
\put(46,15){$\D_{\iq_\is}$}
\put(0,0){first block of¨$\JJ$}
\put(30,0){$\dots$}
\put(55,35){$\dots$}
\put(45,0){last block of¨$\JJ$}
\end{picture}
\caption{\smaller A simple braid divisible by all
$\ss j$'s with $j$ in¨$\JJ$ can be represented by a
diagram finishing with $\D$'s corresponding to the
blocks of¨$\JJ$; the strands going to the same block
cannot cross in the $\bs'$ part, as they cross
inside the block; so the strands starting from¨$i$
and¨$i'$ with $i < i'$ cross if and only if  they
go to different blocks and $i$ goes to a block on
the right of the block $i'$ goes to.}
\label{F:Blocks}
\end{figure}

\begin{proof}
Assume that $\bs$ is a simple $\nn$-braid
satisfying $\DivR(\bs) \supseteq \JJ$, \ie, that is
right divisible by¨$\ss j$ for each¨$j$ in¨$\JJ$.
By hypothesis, the first block of consecutive
elements of¨$\JJ$ is $1, \dots \iq_1-1$. Then $\bs$
being right divisible by $\ss1$, \dots,
$\ss{\iq_1-1}$ is equivalent to its being right
divisible by the left lcm of these elements, which
is¨$\D_{\iq_1}$. Similarly, the second block
in¨$\JJ$ is $\ss{\iq_1+1}$, \dots, $\ss{\iq_1 +
\iq_2-1}$, and being right divisible by these
elements amounts to being right divisible by their
left lcm, which is $\sh^{\iq_1}(\D_{\iq_2})$, where
$\sh$ denotes the shift endomorphism of¨$B_\infty$
that maps each¨$\ss i$ to¨$\ss{i+1}$. Now, by
construction, the blocks $\{1, \dots, \iq_1-1\}$
and $\{\iq_1+1 \dots, \iq_1 + \iq_2-1\}$ are
separated by¨$\iq_1$, and, therefore, the
corresponding¨$\s$'s commute. In particular,
$\D_{\iq_1}$ and $\sh^{iq_1}(\D_{\iq_2})$ commute,
and their left lcm is their product. Finally, a
simple braid¨$\bs$ satisfies
$\DivR(\bs) \supseteq \JJ$ if and only if it is a
right multiple of the element
$$\D_\JJ = \D_{\iq_1} \, \sh^{\iq_1}(\D_{\iq_2}) \,
\dots \, \sh^{\iq_1 + \cdots +
\iq_{\is-1}}(\D_{\iq_\is}),$$ 
\ie, we have $\bs = \bs' \D_\JJ$ for some¨$\bs'$.

We claim that $\ff_\bs$ determines¨$\bs'$,
hence¨$\bs$. Indeed, in a simple braid, any two
strands cross at most once. Now, in a $\D$-diagram,
any two strands cross. So, if the $i$th and the
$i'$th strands go to the same block of¨$\JJ$, \ie,
if we have $\ff_\bs(i) = \ff_\bs(i')$, then these
strands cross in the final¨$\D$-part, and therefore
they cannot cross in (any diagram
representing)¨$\bs'$. So, when
$\ff_\bs$ is given, there is only one way to
construct¨$\bs'$, namely taking the strands to the
entrance of the specified $\D$-block in increasing
order (Figure¨\ref{F:Blocks}). 

Consider now¨$i$, $1 \le i < \nn$. We wonder
whether $\ss i$ is a left divisor of¨$\bs$, \ie, if
the $i$th and the $i + 1$st strands cross in the
diagram of¨$\bs$. If we have $\ff_\bs(i) =
\ff_\bs(i+1)$, the $i$th and $i+1$st strand go to
the same $\JJ$-block, where they certainly
cross. If we have $\ff_\bs(i) >
\ff_\bs(i+1)$, then the $i$th strand goes to a
$\JJ$-block on the right of the $\JJ$-block to
which the $i+1$st strand goes, so they must cross
in the $\bs'$ part. On the contrary, for
$\ff_\bs(i) < \ff_\bs(i+1)$, the strands cannot
cross in the $\bs'$ part---if they crossed once,
they would have to cross a second time before
exiting, and this is forbidden---and they do not
cross in the $\D$ part either. So \eqref{E:Equiv}
holds.
\end{proof}

We deduce the following characterization
of¨$\aas\nn\II\JJ$:

\begin{prop} \label{P:EES}
Assume that $\II, \JJ$ are subsets of~$\Int{\nno}$
with respective $\nn$-composi\-tions $(\ip_1, \dots,
\ip_\ir)$ and $(\iq_1, \dots, \iq_\is)$. Then $\aas\nn\II\JJ$
is the number of $\ir \times \is$ matrices
with nonnegative integer entries such that, for
all~$i, j$, the $i$th row has sum~$\ip_i$ and the
$j$th~column has sum~$\iq_j$. In particular, we have
\begin{equation} \label{E:Value}
\aas\nn\II\emptyset = \frac{\nn!}{\ip_1! \dots
\ip_\ir!}
\mbox{\quad and \quad}
\aas\nn\emptyset\JJ = \frac{\nn!}{\iq_1! \dots
\iq_\is!}.
\end{equation}
\end{prop}

\begin{proof}
Lemma¨\ref{L:Sequence} immediately implies
\begin{gather} 
\label{E:Sequence}
\aa\nn\II\JJ = \card
\{\ff \!:\! \Int\nn \!\to\! \Int\is \,;\, 
(\forall j)(\card\ff\inv(j) = \iq_j) 
\,\&\,
(i \in \II \Leftrightarrow \ff(i) \ge \ff(i+1))\},\\
\label{E:SequenceBis}
\aas\nn\II\JJ = \card\{\ff \!:\! \Int\nn \!\to\!
\Int\is \,;\,  (\forall j)(\card\ff\inv(j) = \iq_j) 
\,\&\,
(i \in \II \Rightarrow \ff(i) \ge \ff(i+1))\}.
\end{gather}

Assume that $\ff$ is a function of~$\Int\nn$
to~$\Int\is$ satisfying the constraints
of~\eqref{E:SequenceBis}. Let $A_\ff$ be the $\ir
\times \is$-matrix whose $(i, j)$-entry is the
number of~$k$'s in the $i$th block of~$\II$
satisfying~$\ff(k) = \iq_j$. By construction, the
the sum of the $i$th row of~$A_\ff$ is the
size~$\ip_i$ of the $i$th block of~$\II$, while
the sum of the $j$th column is the number of~$k$'s
satisfying $\ff(k) = j$, \ie, it is~$\iq_j$. We
claim that $A_\ff$ determines~$\ff$. Indeed,
\eqref{E:SequenceBis} requires that $\ff$ be
non-increasing on each block of~$\II$, so there is
only one possibility once the number of~$k$'s
going to the various~$j$ is fixed.

The first equality in¨\eqref{E:Value} follows:
for¨$\is = \nn$, there is exactly one nonzero entry
in each column, so choosing a convenient matrix
amounts to choosing among $\nn$~elements the
$\iq_1$~columns with a¨$1$ in the first row, the
$\iq_2$~columns with a¨$1$ in the second row, \etc\
The second equality is similar with rows and
columns exchanged. 
\end{proof}

\begin{coro} \label{C:Depends}
$(i)$ The number¨$\aas\nn\II\JJ$ only depends on the
partitions¨$\part\nn\II$ and¨$\part\nn\JJ$.

$(ii)$ For each~$\nn$ and¨$\II$, the number¨$\aa\nn\II\JJ$
only depends on the partition¨$\part\nn\JJ$.
\end{coro}

\begin{proof}
Point¨$(i)$ directly follows from the characterization of
Proposition¨\ref{P:EES}, as the latter clearly involves the
sizes of the blocks of¨$\II$ and¨$\JJ$ only. As for¨$(ii)$, 
the usual inclusion-exclusion formula gives
$$\aa\nn\II\JJ = \sum_{\KK \cap\II = \emptyset}
(-1)^{\card\KK} \aas\nn{\II \cup \KK}\JJ.$$
By¨$(i)$, each term in the sum only depends
on¨$\part\nn\JJ$, and so does the sum.
\end{proof}

It is easy to check that the value of¨$\aa\nn\II\JJ$ does
not only depend on¨$\part\nn\II$ in
general: when we apply the inclusion-exclusion formula, 
the sizes of the blocks in¨$\II \cup \KK$ do not only depend
on the sizes of the block in¨$\II$.

\begin{rema} \label{R:Hohlweg}
Corollary¨\ref{C:Depends} can also be deduced from the classical results by
Solomon about the descent algebra---and, therefore, it extends to all
Artin--Tits groups of spherical type. The argument is as follows. For¨$f$ a
permutation, let
$D(f)$ denote the sets of descents of¨$f$. In the group
algebra¨$\QQ[\Sym\nn]$, let $d_I:= \sum\{f; D(f) = I\}$ and
$e_\JJ= \sum\{f; D(f) \cap \JJ = \emptyset\}$. Using $w_0$
for the flip permutation, we have $w_0 e_\JJ = \sum\{f; D(f) \supseteq
\JJ\}$, and therefore $\aa\nn\II\JJ = \langle d_\II, w_0e_\JJ\rangle$,
where $\langle.,.\rangle$ is the inner product defined by $\langle f,
g\rangle = 1$ for $g = f\inv$, and¨$=0$ otherwise. Using the isometry result
of¨\cite{Hoh}, we deduce $\aa\nn\II\JJ = \langle \theta(\dd_\II),
\theta(w_0e_\JJ)\rangle$, where $\theta(e_\KK)$ denotes the character
of¨$\Sym\nn$ induced by the trivial character of the standard parabolic
subgroup generated by (the transpositions¨$s_i$ with¨$i$ in)¨$\KK$.
By¨\cite{Sol}, the subspace of¨$\QQ[\Sym\nn]$ generated by the¨$e_\KK$ is
a subalgebra, and the kernel of¨$\theta$ is generated by the elements
$\dd_\II - \dd_{\JJ}$ with
$\II, \JJ$ associated with the same partition, and it follows from the above
expression that $\aa\nn\II\JJ$ only depends on the partition
associated with¨$\JJ$.
\end{rema}

\subsection{The matrix¨$\Mattt\nn$}

We can now come back to the matrix¨$\Matt\nn$, and
replace it for the computation of the numbers¨$\bb\nn\dd$
with a new matrix of smaller size. Indeed, a direct
application of Corollary¨\ref{C:Depends} is

\begin{lemm}
Assume that $\JJ$, $\JJ'$ are subsets of¨$\Int\nno$ with
the same $\nn$-partition. Then the $\JJ$th and $\JJ'$th
columns of¨$\Matt\nn$ are equal.
\end{lemm}

Thus the process used to replace¨$\Mat\nn$ with¨$\Matt\nn$
can be applied again, \ie, the form a new matrix by
gathering the equal columns and summing the corresponding
rows.

\begin{defi}
$(i)$ For¨$\l$ a partition (or a composition) of¨$\nn$, we
denote by¨$\set\l$ the unique subset¨$\II$ of¨$\Int\nno$
satisfying $\comp\nn\II = \l$.

$(ii)$ For $\l, \m \vdash \nn$ (\ie, partitions of¨$\nn$),
we put
$$\aab\l\m 
= \sum_{\part\nn\II = \l} \aa\nn\II{\set\m}
= \card\{\bs \,;\, \part\nn{\DivL(\bs)} = \l
\,\&\, \DivR(\bs) \supseteq \set\m\},$$
and we let $\Mattt\nn$ be the matrix with rows and columns
indexed by partitions of¨$\nn$ and whose $(\l, \m)$-entry
is¨$\aab\l\m$.
\end{defi}

In this way the size of the matrix has been reduced
from¨$\nn!$ to¨$\nbpart\nn$, the number of partitions
of¨$\nn$. For instance, enumerating partitions in the order
induced by the previous order on $\Pw(\Int\nn)$, we obtain
$$\Mattt3 = 
\left(\begin{matrix}
1&0&0\\4&2&0\\1&1&1
\end{matrix}\right), \quad
\Mattt4 =
\left(\begin{matrix}1&0&0&0&0\\11&4&1&0&0\\5&3&2&1&0\\
6&4&2&2&0\\1&1&1&1&1
\end{matrix}\right), \quad
\Mattt5 =
\left(\begin{matrix}
1&0&0&0&0&0&0\\
26&8&0&2&0&0&0\\
23&12&4&5&0&1&0\\
43&21&5&10&0&2&0\\
8&6&4&4&2&2&0\\
18&12&6&8&2&4&0\\
1&1&1&1&1&1&1
\end{matrix}\right).$$

Applying the same argument as for
Lemma¨\ref{L:Reduction}, we obtain:

\begin{prop} \label{P:NNNumber}
For $\nn \ge 1$, the characteristic polynomials
of¨$\Mattt\nn$ and $\Mat\nn$ coincide up to a power of¨$\xx$,
and, for every simple $\nn$-braid¨$\bs$, we have for
each¨$\dd$
\begin{equation}
\bbb\nn\dd{\bs} = ((1,1, \dots, 1) \,
\Mattt\nn^{\dd-1})_\l,
\end{equation}
where $\l$ is the $\nn$-partition of¨$\DivL(\bs)$.
In particular we have 
\begin{equation}
\bb\nn\dd = ((1,1, \dots, 1) \,
\Mattt\nn^{\dd})_{(1, 1, \dots, 1)}.
\end{equation}
\end{prop}

Table¨\ref{T:Values} gives the first few values deduced
from the above formulas.

\begin{table}[htb]
$$\begin{tabular}{c|r|r|r|r|r|r}
\quad$\dd$
&1&2&3&4&5&6\\
\hline
\vrule width0pt height12pt depth6pt
$\bbb2\dd1$
& 1
& 2 & 3 & 4 & 5 & 6 \\
\hline
\vrule width0pt height12pt 
$\bbb3\dd1$
& 1
& 6 & 19 & 48 & 109 & 234 \\
\vrule width0pt depth6pt
$\bbb3\dd{\D_2}$
& 1
& 3 & 7 & 15 & 31 & 63 \\
\hline
\vrule width0pt height12pt 
$\bbb4\dd1$
& 1
&24&211&1\,380&8\,077&45\,252
\\
$\bbb4\dd{\D_2}$
& 1
&12&83&492&2\,765&15\,240\\
\vrule width0pt depth6pt
$\bbb4\dd{\D_3}$
& 1
&4&15&64&309&1\,600 \\
\hline
\vrule width0pt height12pt 
$\bbb5\dd1$
& 1
 & 120 & 3\,651 & 79\,140 
& 1\,548\,701 
& 29\,375\,460 
\\
$\bbb5\dd{\D_2}$
& 1
 & 60 & 1\,501 & 30\,540 
& 585\,811
& 11\,044\,080 
\\
$\bbb5\dd{\D_3}$
& 1
 & 20 & 311 & 5\,260 
& 94\,881
& 1\,755\,360 
\\
\vrule width0pt depth6pt
$\bbb4\dd{\D_4}$
& 1
 & 5 & 31 & 325 
& 4\,931 
& 86\,565 
\\
\hline
\vrule width0pt height12pt 
$\bbb6\dd1$
& 1
& 720 
& 90\,921 
& 7\,952\,040 
& 634\,472\,921
& 49\,477\,263\,360
\\
$\bbb6\dd{\D_2}$
& 1
& 360 
& 38\,559 
& 3\,228\,300 
& 254\,718\,389
& 19\,808\,530\,620\\
$\bbb6\dd{\D_3}$
& 1
 & 120 & 8\,727 & 649\,260 
& 49\,654\,757
& 3\,831\,626\,580 \\
$\bbb6\dd{\D_4}$
& 1
 & 30 & 1\,075 & 61\,620 
& 4\,387\,195
& 332\,578\,230 \\
$\bbb6\dd{\D_5}$
& 1
 & 6 & 63 & 1\,955 
& 116\,423
& 8\,448\,606 
\end{tabular}$$
\medskip
\caption{\smaller First values of $\bbb\nn\dd{\D_\rr}$
for $1 \le \rr < \nn$;  the values of¨$\bb\nn\dd$ can be also
read, as Proposition¨\ref{P:Total} gives $\bb\nn\dd =
\bbb\nn{\dd+1}1$, so, for instance, we find $\bb33= 48$, and
$\bb45 = 45\,252$.}
\label{T:Values}
\end{table}

\subsection{Small values of¨$\nn$}

For small values of¨$\nn$, it is easy to complete the
computations and to obtain an explicit form for the
expansion of¨$\bbb\nn\dd\bs$ announced in
Proposition¨\ref{P:Expansion}.

\begin{exam}
Assume $\nn = 3$. The matrix¨$\Mattt3$ is invertible with
eigenvalues¨$1$ (double) and¨$2$. By solving the
recurrences, we find
\begin{equation*} 
\bbb3\dd\bs =
\begin{cases}
4 \cdot 2^\dd - 3\dd -4
&\mbox{¨for $\bs$ with partition $(1, 1, 1)$, \ie,
$\bs = 1$}\\ 
2^\dd - 1, 
&\mbox{¨for $\bs$ with partition $(2, 1)$,
\ie, $\bs = \ss1, \ss2, \ss2\ss1$, or $\ss1\ss2$,}\\
1
&\mbox{¨for $\bs$ with partition $(3)$, \ie, $\bs =
\D_3$,}
\end{cases}
\end{equation*}
and we deduce $\bb3\dd = 8 \cdot 2^\dd - 3\dd - 7$.
\end{exam}

\begin{exam}
Assume now $\nn = 4$. The matrix¨$\Mattt4$ admits 4
eigenvalues, namely those of¨$\Mattt3$, plus $\ev_1 = 3 +
\sqrt6$ and $\ev_2 = 3 - \sqrt6$.  Solving the recurrences
yields for¨$\bbb\nn\dd\bs$ with associated partition as
indicated
$$\begin{array}{lccrccrccrc}
(1,1,1,1):
&\frac1{20}(18 + 7\sqrt6)
&\hspace{-3mm}\ev_1^\dd
&+
&\frac1{20}(18-7\sqrt6)
&\hspace{-3mm}\ev_2^\dd 
&-  
&\frac{256}5 
&\hspace{-3mm}\cdot\, 2^\dd 
&+ 
&6k + 11,\\
(2,1,1):
&\frac1{60}(18 + 7\sqrt6)
&\hspace{-3mm}\ev_1^\dd 
&+
&\frac1{60}(18-7\sqrt6)
&\hspace{-3mm}\ev_2^\dd 
&- 
&\frac85
&\hspace{-3mm}\cdot\, 2^\dd 
&+ 
&1,\\
(2,2):
&\frac1{12}(\sqrt6)
&\hspace{-3mm}\ev_1^\dd  
&- 
&\frac1{12}(\sqrt6)
&\hspace{-3mm}\ev_2^\dd,\\
(3,1):
&\frac1{60}(6 -\sqrt6)
&\hspace{-3mm}\ev_1^\dd
&+
&\frac1{60}(6+\sqrt6)
&\hspace{-3mm}\ev_2^\dd 
&+ 
&\frac45
&\hspace{-3mm}\cdot\, 2^\dd 
&-
&1,\end{array}
$$
and¨$1$ for associated partition $(4)$, \ie, for $\bs =
\D_4$. As the characteristic polynomial of¨$\Mattt4$ is
$(x^2 - 6x + 3)(x - 2)(x-1)^2$,
we can equivalently determine $\bbb4\dd\bs$ and $\bb4\dd$ by
inductions on¨$\dd$ of the form
\begin{equation} \label{E:Recur}
\var_\dd = 6\var_{\dd-1} - 3\var_{\dd-2}
+ \alpha 2^\dd + \beta \dd + \gamma
\end{equation}
where $\alpha, \beta, \gamma$ are determined using special
values of¨$\var_\dd$. For instance, $\bb4\dd$ is determined
by¨\eqref{E:Recur} with $\alpha = 32$, $\beta = -12$,
$\gamma = -34$ and the values $\var_{-1} = 0$, $\var_0= 1$. 
Generating functions can be deduced easily.
\end{exam}

\subsection{Eigenvalues of $\Mat\nn$}

By Proposition¨\ref{P:Expansion}, the value of¨$\bb\nn\dd$
and¨$\bbb\nn\dd\bs$, and in particular its asymptotic
behaviour when $\dd$ grows to infinity, are connected with
the non-zero eigenvalues of¨$\Mat\nn$, which, by
Proposition¨\ref{P:NNNumber}, coincide with those
of¨$\Mattt\nn$. The characteristic polynomial
of¨$\Mattt\nn$---hence of¨$\Mat\nn$ up to an
$x^\dd$¨factor---for small values of¨$\nn$ is displayed in
Table¨\ref{T:Charpol}.

\begin{table}[t]
\begin{tabbing}
\hspace{10mm}\=\hspace{85mm}\=\hspace{2cm}\kill
\> $\charpol{\Mattt1} 
= x - 1$\\
\> $\charpol{\Mattt2} 
= \charpol{\Mattt1}  \cdot  (x - 1)$\\
\> $\charpol{\Mattt3} 
= \charpol{\Mattt2}  \cdot  (x - 2)$\\
\> $\charpol{\Mattt4} 
= \charpol{\Mattt3}  \cdot  (x^2 - 6x +3)$\\
\> $\charpol{\Mattt5} 
=\charpol{\Mattt4}  \cdot  (x^2 - 20 x +24)$\\
\> $\charpol{\Mattt6} 
= \charpol{\Mattt5}
 \cdot  (x^4 - 82 x^3 +359x^2 - 260x + 60)$\\
\> $\charpol{\Mattt7} 
= \charpol{\Mattt6} \cdot
(x^4  - 390 x^3  + 6024 x^2 - 13680 x + 8640)$\\
\> $\charpol{\Mattt8} 
= \charpol{\Mattt7}
\cdot (x^7  - 2134 x^6  + 139976 x^5 - 1321214 x^4
+ 3780975 x^3$\\
\> \hspace{6cm}$- 3305160 x^2  + 1341900 x - 226800)$
\end{tabbing}
\begin{tabular}{c|c|c|c|c|c|c|c|c}
$\nn$&1&2&3&4&5&6&7&8\\
\hline
\vline width0pt height10pt depth6pt
$\ev_{max}(\Mat\nn)$&1&1&2&5.449 &18.717 
&77.405 &373.990 &2066.575 \\
$\displaystyle\frac{\ev_{max}(\Mat\nn)}{\nn\cdot\ev_{max}(\Mat\nno)}$
&-&0.5& 0.667 &0.681 &0.687 &0.689 &0.690 &0.691\\ 
\end{tabular}
\bigskip
\caption{\smaller Characteristic polynomial
of¨$\Mattt\nn$ for $\nn \le 8$, and the
corresponding largest eigenvalue---which is to be
compared with¨$\nn!$, the growth rate
for the number of $\nn$-braids of degree at
most¨$\dd$ if all sequences were normal.}
\label{T:Charpol}
\end{table} 

These values support the following

\begin{conj} \label{C:Inclusion}
For each¨$\nn$, the characteristic polynomial
of¨$\Mat\nno$ divides that of¨$\Mat\nn$. More precisely,
the sprectrum of¨$\Mattt\nn$ is the spectrum
of¨$\Mattt\nno$, plus $p(\nn) - p(\nno)$ simple non-zero
eigenvalues.
\end{conj}

It is not hard to check the above statement for $\nn \le
10$. Specifically, let¨$\widehat M_\nn$ be the size
$p(\nn)-2$ matrix obtained from¨$\Mattt\nn$ by deleting the
first and the last rows and columns. For each small value
of¨$\nn$, one can directly check that $\widehat M_\nn$ is
similar to a matrix of the form $\left(\begin{matrix}
\widehat{M}_\nno&0\\\dots&\dots\end{matrix}\right)$, and
deduce the properties asserted in
Conjecture¨\ref{C:Inclusion}. But ne generic argument is
known so far.

The growth rate of the numbers¨$\bbb\nn\dd\bs$
is connected with the largest eigenvalue
$\ev_{max}(\Mat\nn)$ of¨$\Mat\nn$. For¨$\nn \le 6$,
all¨$\bbb\nn\dd\bs$ except¨$\bbb\nn\dd{\D_\nn}$,
which is¨$1$, and therefore all $\bb\nn\dd$ as
well, grow like $\ev_{max}(\Mat\nn)^\dd$. 

\begin{ques}
Do all¨$\bbb\nn\dd\bs$ except¨$\bbb\nn\dd{\D_\nn}$
grow like $\ev_{max}(\Mat\nn)^\dd$?
\end{ques}

\begin{ques}
What is the asymptotic behaviour of $\ev_{max}(\Mat\nn)$
with¨$\nn$?
\end{ques}

The trivial upper bound of¨\eqref{E:Upper} suggests to
compare $\ev_{max}(\Mat\nn)$ with¨$\nn!$, or, rather,
$\ev_{max}(\Mat\nn)$ with¨$\nn \cdot \ev_{max}(\Mat\nno)$.

\section{Letting the degree vary}
\label{S:Deg}

So far, we kept the braid index¨$\nn$ fixed, and studied
how the numbers¨$\bb\nn\dd$ or¨$\bbb\nn\dd\bs$ vary
with¨$\dd$, thus letting linear inductions appear. Quite
different induction schemes appear when we fix the degree
and let the braid index vary. No systematic method is
known so far, and we only mention a few partial results
motivated by the approach of¨\cite{Dhh}.

\subsection{The numbers¨$\bb\nn2$}

Very little is known about¨$\bb\nn\dd$ in general. The case
$\dd = 1$ is trivial, as we already observed the equality 
$$\bb\nn1 = \nn!.$$
For $\dd = 2$, the value can be deduced from earlier results
of¨\cite{CSV, CSV2} about permutations. We shall use the
following very general observation about duality in Garside
groups:

\begin{lemm} \label{L:Dual}
For¨$\xx$ in¨$\DD\nn{}$, let $\dual\xx$
and $\xx^*$ be defined by $\dual\xx\,\xx = \xx \,
\xx^* = \D_\nn$. Then $\xx \mapsto \dual\xx$ and
$\xx \mapsto \xx^*$ are permutations of¨$\DD\nn{}$,
and, for each simple¨$\xx$, we have
\begin{equation} \label{E:Dual}
\DivR(\dual\xx) = \Int\nn \setminus \DivL(\xx)
\mbox{\quad and \quad}
\DivL(\xx^*) = \Int\nn \setminus \DivR(\xx).
\end{equation}
\end{lemm}

\begin{proof}
Assume $\xx \in \DD\nn{}$. Then, by hypothesis,
$\xx$ is a left and a right divisor of¨$\D_\nn$,
hence $\dual\xx$ and $\xx^*$ are positive
braids, and they are divisors of¨$\D_\nn$
in¨$B_\nn^+$, so they are simple. That the
mappings $\xx \mapsto \dual\xx$ and
$\xx \mapsto \xx^*$ are injective is clear,
and the surjectivity follows from the
finiteness of¨$\DD\nn{}$.

Now, $\ss i$ being a right divisor of¨$\dual\xx$ is
equivalent to $\dual\xx \ss i$ not being simple,
hence to the non-existence of¨$\yy$ satisfying
$\dual\xx \ss i \yy = \D_\nn$, and finally to the
non-existence of¨$\yy$ satisfying $\xx = \ss i
\yy$. This implies the first equality
in¨\eqref{E:Dual}. The second equality follows from
a symmetric argument.
\end{proof}

\begin{prop}
The numbers¨$\bb\nn2$ are determined by the
induction
\begin{equation}
\bb02=1, \qquad 
\bb\nn2 = \sum_{\ind=0}^\nno (-)^{\nn + \ind
+1} {\nn \choose \ind}^2 
\bb\ind2.
\end{equation}
Their double exponential generating function is 
\begin{equation}
\sum_{\nn = 0}^\infty \bb\nn2
\frac{x^\nn}{\nn!^2} = 
\bigg( \sum_{\nn = 0}^\infty (-1)^\nn
\frac{x^\nn}{\nn!^2}
\bigg)\inv = \frac1{J_0(\sqrt\xx)},
\end{equation}
where $J_0(\xx)$ is the Bessel function.
\end{prop}

\begin{proof}
By definition, $\bb\nn2$ is the number of pairs
of simple¨$\nn$-braids $(\sx_1, \sx_2)$ satisfying
$\DivR(\sx_1) \supseteq \DivL(\sx_2)$, \ie, by
Lemma¨\ref{L:Dual},
$\DivR(\sx_1) \cap \DivR(\dual\sx_2) = \emptyset$. By
Lemma¨\ref{L:Exchange}, this number is also
the number of pairs of permutations¨$(\ff, \gg)$
in¨$\Sym\nn$ with no descent in common, \ie, such
that there exists no¨$i$ satisfying both $\ff(i) >
\ff(i+1)$ and $\gg(i) > \gg(i+1)$. Such pairs of
permutations have been counted
in¨\cite{CSV, CSV2} (see also¨\cite{Rio}),
with the result indicated above.
\end{proof}

\subsection{The numbers $\bbb\nn2{\D_{\nn - \rrr}}$}

Specific results appear when we consider the
numbers¨$\bbb\nn\dd\bs$ with¨$\bs$ of the form¨$\D_{\nn
- \rrr}$ with $1 \le \rrr \le \nn$. In particular, we can
complete the computation when $\rrr$ is fixed and $\dd$ is
small. We obviously have
$\bbb\nn1{\D_{\nn - \rrr}} = 1$ for $1 \le \rrr \le \nn$, so
the first case to consider is $\dd = 2$. The general
principle that makes the computation of $\bbb\nn\dd{\D_{\nn
- \rrr}}$ relatively easy is the following observation:

\begin{lemm} \label{L:Prelim}
For all¨$\nn, \dd, \rrr$, we have
\begin{equation} \label{E:Prelim}
\bbb\nn\dd{\D_{\nn - \rrr}} = \sum_{\mbox{\Small\rm $\bs$ right divisible
by¨$\D_{\nn - \rrr}$}} \bbb\nn{\dd-1}\bs.
\end{equation}
\end{lemm}

\begin{proof}
The argument is similar to that for
Proposition¨\ref{P:Total}. A sequence $(\bs_1,
\dots, \bs_{\dd-1}, \D_{\nn-\rr})$ is normal if and only if
both $(\bs_1, \dots, \bs_{\dd-1})$ and $(\bs_{\dd-1},
\D_{\nn-\rr})$ are normal. Now $(\bs_{\dd-1}, \D_{\nn-\rr})$
is normal if and only if every¨$\ss i$
dividing¨$\D_{\nn-\rr}$ on the left divides¨$\bs_{\dd-1}$
on the right. The $\ss i$'s dividing $\D_{\nn-\rr}$ on the
left are $\ss1$, \dots, $\ss{\nn - \rr - 1}$. The simple
braids that are right divisible by $\ss1$, \dots, $\ss{\nn -
\rr - 1}$ are those right divisible by¨$\D_{\nn-\rr}$. Then
\eqref{E:Prelim} follows.
\end{proof}

\begin{prop}
For $1 \le \rrr \le \nn$, we have
\begin{equation} \label{E:K1}
\bbb\nn2{\D_{\nn - \rrr}} = \frac{\nn!}{(\nn - \rrr)!}.
\end{equation}
\end{prop}

\begin{proof}
By Lemma¨\ref{L:Prelim}, $\bbb\nn2{\D_{\nn - \rrr}}$ is
the number of simple $\nn$-braids¨$\bs$ that are
right divisible by¨$\D_{\nn - \rrr}$, \ie., that satisfy
$\DivR(\bs) \supseteq \Int{{\nn - \rrr}}$. The block
composition of¨$\Int{\nn - \rrr}$ in¨$\nn$ is $({\nn - \rrr}, 1,
\dots, 1)$, so \eqref{E:Value} directly
gives¨\eqref{E:K1}.
\end{proof}

\subsection{The numbers $\bbb\nn3{\D_{\nn -\rrr}}$}

Things become more interesting for¨$\dd = 3$. 

\begin{prop} \label{P:K2}
For $1 \le \rrr \le \nn$, there exist polynomials¨$P_1$, \dots, $P_{n-r}$
with integer coefficients and $P_i$ of degree at most¨${\nn - \rrr} - i +
1$ such that, for every¨$\nn$, we have
\begin{align} 
\label{E:K2}
\bbb\nn3{\D_{\nn-\rrr}} 
&= (\nn - \rrr)!\,({\nn - \rrr}+1)^\nn + \sum_{i = 1}^{\nn -
\rrr} P_i(\nn) \, i^{\rrr +i-1}.\\
\intertext{The explicit values for $\rrr
= 1, 2$ are}
\label{E:K21}
\bbb\nn3{\D_{\nno}} &= 2^\nno, \\
\label{E:K22}
\bbb\nn3{\D_{\nn-2}} &= 2 \cdot 3^\nn - (\nn + 6) \cdot
2^\nno + 1.
\end{align}
\end{prop}

\begin{proof}
We begin with¨\eqref{E:K21}. By
Lemma¨\ref{L:Prelim}, $\bbb\nn3{\D_{\nno}}$ is the sum of
all $\bbb\nn2\bs$ with $\bs$ right divisible by¨$\D_{\nno}$, \ie, it is the number of normal sequences
$(\sx_1, \sx_2)$ such that $\sx_2$ is right divisible
by¨$\D_{\nno}$. Let $\SS$ be the set of all such
normal sequences. We partition¨$\SS$ according to the
value of¨$\DivL(\sx_2)$, \ie, for each subset¨$\II$
of¨$\Int{\nno}$, we count how many pairs¨$(\sx_1,
\sx_2)$ satisfy¨$\DivL(\sx_2) = \II$. So assume
that $\sx_2$ is right divisible by¨$\D_{\nno}$. Two
cases are possible. Either $\sx_2$ is right
divisible by (hence equal to)¨$\D_\nn$, and then we
have
$\DivL(\sx_2) = \Int{\nno}$. Or
$\sx_2$ is not divisible by¨$\D_{\nno}$, and then
Lemma¨\ref{L:Sequence} shows that
$\sx_2$ must be $\ss{i, \nn}\D_{\nno}$ for some¨$i$
with $2 \le i \le \nn$, so that the
block composition of¨$\DivL(\sx_2)$ is $(i-1,
\nn - i+1)$. So the possible compositions
for the set¨$\DivL(\sx_2)$ are¨$(\nn)$, and $(\ip,
\nn - \ip)$ with $1 \le \ip \le \nno$. Conversely,
the previous analysis shows that, for each¨$\II$ of
the previous form, there exists exactly one
possible¨$\sx_2$. Now, Proposition¨\ref{P:EES} says
that there is one choice for¨$\sx_1$ in the case
of¨$(\nn)$---namely $\sx_1 = \D_\nn$---and $\nn
\choose \ip$ choices for¨$\sx_1$ in the case
of¨$(\ip, \nn - \ip)$. We deduce
$$\bbb\nn3{\D_{\nno}} = 1 + \sum_{\ip= 1}^\nno {\nn
\choose{\ip}} = 2^\nno.$$

The method is similar for
computing¨$\bbb\nn3{\D_{\nn-2}}$ in¨\eqref{E:K22}. Assume
that
$(\sx_1, \sx_2)$ is a normal sequence with¨$\sx_2$ right
divisible by¨$\D_{\nn-2}$. The hypothesis is $\DivR(\sx_2) \supseteq
\Int{\nn-3}$, so three cases may occur, namely $\DivR(\sx_2)
\supseteq \Int{\nn-2}$, $\DivR(\sx_2) = \Int{\nn-3}$, and
$\DivR(\sx_2) =  \Int{\nn-3} \cup \{\nno\}$. The first case was analysed
above. In the second case, $\DivL(\sx_2)$ has three blocks, and,
conversely, each set¨$\II$ with three blocks gives
exactly one eligible¨$\sx_2$. In the third case,
$\DivL(\sx_2)$ has either two blocks, or it has
three blocks with the middle one at least¨$2$;
conversely, each set¨$\II$ of the previous form
gives one eligible¨$\sx_2$. Using as above
Proposition¨\ref{P:EES} to count the
eligible¨$\sx_1$'s for each possible¨$\II$, we
obtain that $\bbb\nn3{\D_{\nn-2}}$ is
\begin{equation} \label{E:Sum}
\bbb\nn3{\D_{\nno}} 
+ \sum_{\overset{\scriptstyle\ip_1 + \ip_2 + \ip_3 =
\nn} {\ip_1, \ip_2, \ip_3 \ge 1}}
\frac{\nn!}{\ip_1! \ip_2! \ip_3!}
+ \sum_{\overset{\scriptstyle\ip_1 + \ip_2 =
\nn} {\ip_1, \ip_2 \ge 1}}
\frac{\nn!}{\ip_1! \ip_2!}
+ \sum_{\overset{\scriptstyle\ip_1 + \ip_2 + \ip_3 =
\nn} {\ip_1, \ip_3 \ge 1,\ip_2 \ge 2}}
\frac{\nn!}{\ip_1! \ip_2! \ip_3!}.
\end{equation}
Using the fact that $3^\nn$ is the sum of all
$\frac{\nn!}{\ip_1! \ip_2! \ip_3!}$ with $\ip_1
+ \ip_2 + \ip_3 = \nn$, one deduces¨\eqref{E:K22}
by bookkeeping.

Applying the same method in the general case leads
to¨\eqref{E:K2}. Indeed, always by
Lemma¨\ref{L:Sequence}, specifying a simple
$\nn$-braid¨$\sx_2$ satisfying¨$\DivR(\sx_2)
\supseteq \Int{\rrr -1}$ amounts to choosing a
permutation of the ${\nn - \rrr}$¨last strands and the
${\nn - \rrr}$¨positions $(i_1, \dots, i_{\nn - \rrr})$ where these
strands start from. In the generic case, the
resulting set¨$\DivL(\sx_2)$ is $\{i_1 - 1, \dots,
i_{\nn - \rrr} - 1\}$, whose composition consists of
${\nn - \rrr}$¨blocks. The special cases are when at least
two adjacent strands among the last ${\nn - \rrr}$¨ones
start from adjacent positions; according to
whether these strands cross or not in the final
part, one then obtains either a composition with a
block of size¨2 at least, or a composition with
less than ${\nn - \rrr}$¨blocks. Conversely, for every
subset¨$\II$ of¨$\Int{\nno}$ with ${\nn - \rrr}$¨blocks,
there exists in general $(\nn - \rrr)!$¨eligible¨$\sx_2$'s,
one for each choice of the final permutation of the
last $\nn - \rrr$¨strands. There may be less than
$(\nn - \rrr)!$¨choices for¨$\sx_2$ when $1$ occurs in the
composition of¨$\II$. Also, subsets
of¨$\Int{\nno}$ with fewer than¨${\nn - \rrr}$ blocks
may lead to eligible¨$\sx_2$'s. Multiplying by the
number of eligibles¨$\sx_1$'s for each¨$\II$ and
summing up yields an expression similar
to¨\eqref{E:Sum}, involving $(\nn - \rrr)!$¨sums of the form
$\sum_{\ip_1 + \dots + \ip_{{\nn - \rrr}+1} = \nn}
\frac{\nn!}{\ip_1! \dots \ip_{{\nn - \rrr}+1}!}$ with
possible order constraints on¨$\ip_1$, \dots,
$\ip_{{\nn - \rrr}+1}$. Each of them leads to a
factor¨$({\nn - \rrr}+1)^\nn$, plus additional factors
corresponding to specializing arguments to¨$0$
or¨$1$ or to grouping them.
\end{proof}

\subsection{The numbers $\bbb\nn4{\D_{\nno}}$}

For $\dd = 4$, it seems hopeless to complete the computation
of¨$\bbb\nn\dd{\D_{\nn-\rr}}$. However, this can be done
for¨$\rrr = 1$. The remarkable point is that still another
induction scheme appears.

\begin{prop}
For¨$\nn \ge 1$, we have
\begin{equation} \label{E:KKK3}
\bbb\nn4{\DD{\nno}{}} = \sum_{\ind=0}^\nno
\frac{\nn!}{\ind!}.
\end{equation}
\end{prop}

\begin{proof}
According to Lemma¨\ref{L:Prelim} again, we have
now to count the normal sequences
$(\sx_1, \sx_2, \sx_3)$ with¨$\sx_3$ of the
form¨$\ss{i, \nn}\D_{\nno}$, $2 \le i \le \nn$.
We partition the family according to the value¨$\II$
of¨$\DivL(\sx_2)$, and count how many sequences may
correspond to a given¨$\II$. Let $(\ip_1, \dots,
\ip_\ir)$ denote the block composition of¨$\II$.

Let us first consider the case $\II =
\Int{\nno}$. Then we must have $\sx_2 = \D_\nn$,
hence $\sx_1 = \D_\nn$ as well. There are
$\nn$¨possible choices for¨$\sx_3$, and the total
number of corresponding sequences¨$(\sx_1, \sx_2,
\sx_3)$ is¨$\nn$.

We assume now $\II \not= \Int{\nno}$, \ie, $\ir
\ge 2$. As for¨$\sx_1$, Proposition¨\ref{P:EES} directly
gives the number of choices, namely
$\frac{\nn!}{\ip_1! \dots \ip_\ir!}$. So we are left
with counting how many pairs¨$(\sx_2, \sx_3)$ are
eligible. The case $\sx_3 = \D_n$ is excluded
since it implies $\sx_2 = \D_\nn$ hence $\II =
\Int{\nno}$. As in the case of¨$\bbb\nn3{\D_{\nno}}$,
the hypothesis that $\sx_3$ is $\ss{i, \nn}
\D_{\nno}$ for some¨$i$ with $2 \le i \le \nn$
implies that the block composition
of¨$\DivL(\sx_3)$ consists of two nonempty blocks,
and, conversely, each partition of¨$\Int\nn$ into
two nonempty blocks gives a unique¨$\sx_3$ of the
convenient form. So the number of pairs¨$(\sx_2,
\sx_3)$ associated with¨$\II$ is the number
of¨$\sx_2$'s satisfying
$\DivL(\sx_2) = \II$ and such that $\DivR(\sx_2)$
has two blocks.

By¨\eqref{E:Sequence}, this number is the
number of functions¨$\ff$ of¨$\Int\nn$ to¨$\{1,
2\}$ such that $\ff(i) < \ff(i+1)$ holds exactly for
$i \notin \II$. As only two values are possible,
this condition means that we have $\ff(i) = 1$ and
$\ff(i+1) = 2$ for¨$i \notin \II$, and
$\ff(i+1) \le \ff(i)$ for¨$i \in \II$. Consider
the blocks of¨$\II$. In each block,
except possibly the first and the last ones, the
value of¨$\ff$ has to be¨$2$ on the first element,
and to be¨$1$ on the last element. Inbetween, $\ff$
is non-increasing. So the values consist of a series
of¨$2$'s, followed by a series of¨$1$'s. The
only parameter to specify is the position
where¨$\ff$ switches from¨$2$ to¨$1$, so, for a
block of size¨$\ip$, there are $\ip-1$ possible
choices (see Figure¨\ref{F:Function}). The cases of
the first and the last blocks are special, because
there is no constraint on the left for the first
block, and on the right for the last block. So, in
these special cases, there are
$\ip$¨choices instead of¨$\ip-1$. The conclusion is
that, for¨$\II$ of block composition $(\ip_1, \dots,
\ip_\ir)$, there are
$\ip_1(\ip_2-1) \dots (\ip_{\ir-1}-1) \ip_\ir$
choices for the pairs¨$(\sx_2, \sx_3)$ associated
with¨$\II$. Merging the result for¨$\sx_1$ and
for¨$(\sx_2, \sx_3)$ and summing up over¨$\II$
gives
\begin{equation} \label{E:K3}
\bbb\nn4{\DD{\nno}{}} = \sum
\frac{\nn!}{\ip_1! \dots \ip_\ir!} \,  
\ip_1(\ip_2-1)
\dots (\ip_{\ir-1}-1)\ip_\ir.
\end{equation}
the sum being taken over all finite compositions
$(\ip_1, \dots, \ip_\ir)$ of¨$\nn$: indeed, the
value for¨$\II = \Int{\nno}$, namely¨$\nn$,
corresponds to the missing term $\frac{\nn!}{\nn!}
\nn$ of the sum.

\begin{figure} [htb]
\begin{picture}(82,20)(0, 0)
\put(0,3){\includegraphics{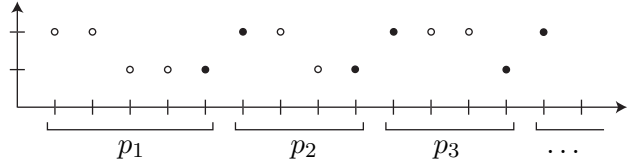}}
\put(15,0){$\ip_1$}
\put(38,0){$\ip_2$}
\put(57,0){$\ip_3$}
\put(72,0){$\dots$}
\end{picture}
\caption{\smaller Proof of¨\eqref{E:K3}: the
rises are fixed, so it just remains to choose the
position of the fall in each block of¨$\II$,
whence $\ip-1$ choices for a size¨$\ip$ block
except the first and the last ones.}
\label{F:Function}
\end{figure}

We can now simplify the right hand term
in¨\eqref{E:K3}. To this end, we observe that
\begin{equation} \label{E:KK3}
\sum_{\underset{\scriptstyle \ip_1, \dots, \ip_\ir \ge 1}{\ip_1 +
\cdots + \ip_\ir = \ind + 1}}
\frac{\ip_1-1}{\ip_1!} \dots \frac{p_{r-1}-1}{p_{r-1}!} \,
\frac{\ip_\ir}{\ip_\ir!} = 1
\end{equation} 
holds for¨$\ind \ge 0$. Indeed, let $F(\ind)$ be the
left hand side of¨\eqref{E:KK3}. We
prove¨\eqref{E:KK3} using induction on¨$\ind$. For
$\ind=0$, we get $1=1$. Assume $\ind \ge 1$ and
consider the sequences
$(\ip_1, \dots, \ip_\ir)$ satisfying $\ip_1 + \cdots + 
\ip_\ir = \ind + 1$.  On the one hand, we have
$(\ind + 1)$, whose contribution to¨$F(\ind)$ is
$\frac{1}{\ind!}$. On the other hand, we have
the sequences of length at least¨$2$. Now, for
each¨$\ip$ with $1 \le \ip \le \ind$, the
contribution of $(\ip, \ip_2, \dots, \ip_\ir)$
to¨$F(\ind)$ is
$\frac{\ip-1}{\ip!}$ times the contribution
of $(\ip_2, \dots, \ip_\ir)$ to¨$F(\ind - \ip)$.
Hence the total contribution of the sequences
beginning with¨$\ip$ to¨$F(\ind)$ is
$\frac{\ip-1}{\ip!}
\, F(\ind - \ip)$, so, by induction hypothesis, it
is $\frac{(\ip-1)}{\ip!}$. We deduce
$F(\ind) = \frac{0}{1!} + \frac{1}{2!} + \cdots +
\frac{\ind -1}{\ind!} + \frac1{\ind!}$,  
which is clearly¨$1$. 

Consider now the right hand side in¨\eqref{E:K3}.
For
$0 \le \ind < \nno$, the contribution of $(\ind+1,
\ip_2, \dots, \ip_\ir)$ to the sum is
$\frac{\nn!}{\ind!}$ times the  quantity
$\frac{\ip_1-1}{\ip_1!} \,
\frac{\ip_{\ir_1}-1}{\ip_{\ir-1}!}
\, \frac{\ip_\ir}{\ip_\ir!}$ involved
in¨\eqref{E:KK3}. Using the latter equality, we
deduce that the total contribution of the sequences
beginning with¨$\ind + 1$ is
$\frac{\nn!}{\ind!}$. As for¨$\ind = \nno$, the
contribution of¨$(\nn)$ to the right
hand side in¨\eqref{E:K3} is¨$\nn$, which is
$\frac{\nn!}{(\nno)!}$, so the general formula
remains valid. By summing over¨$\ind$, we
obtain¨\eqref{E:KKK3}.
\end{proof}

\begin{coro}
The numbers $\bbb\nn4{\DD{\nno}{}}$ are
determined by the induction
$$\var_1 = 1, \qquad
\var_\nn = \nn \var_{\nno} +
2\nno.$$
\end{coro}

Another consequence of¨\eqref{E:KKK3} is the equality
$$\bbb\nn4{\DD{\nno}{}} = \lfloor \nn! e \rfloor - 1,$$
with $e = \exp(1)$.

\end{document}